\documentstyle[11pt]{article}
\newcommand{\B}{{\cal B}}

\newcommand{\C}{{\cal C}}

\newcommand{\A}{{\cal A}}

\newcommand{\G}{{\cal G}}

\newtheorem{theorem}{Theorem}[section]

\newtheorem{lemma}[theorem]{Lemma}

\newtheorem{construction}[theorem]{Construction}
\def\whitebox{{\hbox{\hskip 1pt
 \vrule height 6pt depth 1.5pt
 \lower 1.5pt\vbox to 7.5pt{\hrule width
    3.2pt\vfill\hrule width 3.2pt}%
 \vrule height 6pt depth 1.5pt
 \hskip 1pt } }}
\def\qed{\ifhmode\allowbreak\else\nobreak\fi\hfill\quad\nobreak
     \whitebox\medbreak}
\newcommand{\proof}{\noindent{\it Proof:}\ }

\newcommand{\ignore}[1]{}

\begin{document}
\title{The completion of optimal $(3,4)$-packings\thanks{Research supported by NSFC grant 11222113.}}

 \author{{\small Jingjun Bao and Lijun Ji\thanks{Corresponding author}}\\
 {\small Department of Mathematics, Soochow University, Suzhou
 215006, China}\\
 {\small E-mail: jilijun@suda.edu.cn}\\
 }

\date{}
\maketitle
\begin{abstract}
\noindent \\
 A 3-$(n,4,1)$ packing design consists of an
$n$-element set $X$ and a collection of $4$-element subsets of
$X$, called {\it blocks}, such that every $3$-element subset of
$X$ is contained in at most one block. The packing number of
quadruples $d(3,4,n)$ denotes the number of blocks in a maximum
$3$-$(n,4,1)$ packing design, which is also the maximum number
$A(n,4,4)$ of codewords in a code of length $n$, constant weight
$4$, and minimum Hamming distance 4.  In this paper the undecided 21
packing numbers $A(n,4,4)$ are shown to be
equal to Johnson bound $J(n,4,4)$ $(
=\lfloor\frac{n}{4}\lfloor\frac{n-1}{3}\lfloor\frac{n-2}{2}\rfloor\rfloor\rfloor)$
 where $n=6k+5$,
 $k\in \{m:\ m$ is odd, $3\leq m\leq
35,\ m\neq 17,21\}\cup \{45,47,75,77,79,159\}$.

\medskip


\end{abstract}

\section{Introduction}

A 3-$(n,4,1)$ {\it packing design} consists of an $n$-element set
$X$ and a collection of $4$-element subsets of $X$, called {\it
blocks}, such that every $3$-element subset of $X$ is contained in
at most one of them. Such a design is called a {\it packing
quadruple} and denoted by PQS$(n)$ (as in \cite{HP1992}).

A PQS$(n)$ $(X,{\cal A})$ is called {\it maximum} if there does
not exist any PQS$(n)$ $(X,{\cal B})$ with $|{\cal A}|<|{\cal
B}|$, and shortly denoted by MPQS$(n)$. The packing number is the
number of blocks in an MPQS$(n)$ and denoted by $d(3,4,n)$, and by
$A(n,4,4)$, where $A(n,d,w)$ is the maximum number of codewords in
a code of length $n$, constant weight $w$, and minimum Hamming
distance $d$.

 The problem of determining $A(n,4,4)$ has received a lot of
attention from the point of view of combinatorics and coding
theory.

It is known that the Johnson bound $J(n,4,4)$ for the packing
numbers \cite{Johnson} is given by

\[
A(n,4,4)\leq J(n,4,4)=\left\{
\begin{array}{ll}
\lfloor\frac{n}{4}\lfloor\frac{n-1}{3}\lfloor\frac{n-2}{2}\rfloor\rfloor\rfloor & n\not\equiv 0\pmod 6,\medskip \\
\lfloor\frac{n}{4}\lfloor\frac{n-1}{3}\lfloor\frac{n-2}{2}\rfloor\rfloor-1\rfloor
& n\equiv 0\pmod 6.
\end{array}
\right.
\]

\noindent Here, $\lfloor x\rfloor$ denotes the largest integer not
more than $x$.

When $n\equiv 2,4\ (mod\ 6)$, Hanani \cite {Hanani1960} showed
that $A(n,4,4)=J(n,4,4)$ by constructing a PQS$(n)$ with the
property that each triple is contained in exactly one block. Such
a design is called a {\it Steiner quadruple system} of order $n$
and denoted by SQS$(n)$. Deleting one point and all blocks
containing it from an SQS$(n+1)$ yields that $A(n,4,4)=J(n,4,4)$
if $n\equiv 1,3\pmod 6$. Brouwer \cite {Brouwer1978} showed
$A(n,4,4)=J(n,4,4)$ for $n\equiv 0\pmod 6$. The second author showed that $A(n,4,4)=J(n,4,4)$ for $n\equiv 5\pmod 6$ with 21 possible values \cite{Ji2006}.
These results are
summarized as follows.

\begin{theorem}{\rm \cite{Brouwer1978,Hanani1960,Ji2006}}
\label{1.1} For any positive integer $n \not\in  \{6k+5:k=3,5,7,9,11,13,15,19,23,25,$ $27,29,31,33,35,45,47,75,77,79,159\}$, $A(n,4,4)=J(n,4,4)$.
\end{theorem}

The purpose of this paper is to determine the last 21 undecided packing numbers
$A(n,4,4)$. Throughout the remainder of this paper, an MPQS$(n)$ is always
assumed to have $J(n,4,4)$ blocks.

The rest of this paper is arranged as follows. In Section 2, we
construct an MPQS$(n)$ for $n\in \{23,35,47,59,71\}$ directly. In Section
3, we describe recursive constructions for MPQS$(n)$'s via candelabra
quadruple systems. In Section 4 we determine  the last 21 undecided packing numbers
$A(n,4,4)$. Combining these results with Theorem \ref{1.1}, the packing numbers
$A(n,4,4)$ are then completely determined.

\section{Small values}

In this section we  construct an MPQS$(n)$ for $n\in \{23,35,47,59,71\}$.

\begin{lemma}
\label{MPQS(23)} There is an MPQS$(23)$.
\end{lemma}

\proof Let $X=\{0,1,2,\ldots,22\}$ and let $\alpha$ be a permutation as follows.
$$\alpha=(0\ 1)(2\ 3\ 4)(5\ 6\ 7\ 8\ 9\ 10)(11\ 12\ 13\ 14\ 15\ 16)(17\ 18\ 19\ 20\ 21\ 22)$$
The following base blocks generate the required $J(23,4,4)=419$ blocks under the action of the permutation $\alpha$, where the first one base block generates only two distinct blocks and each of the other five base blocks in the first row generates three distinct blocks.

  \begin{center}
\begin{tabular}{llllllllllllll}
  0  5  7  9  & 0  1  5  8  &  0  1 11 14 &   0  1 17 20 &   3  4  5  8 &   5  6  8  9   \\
  0  2  3  6  &  0  2  5 10 &   0  2  7 11 &   0  2  9 12 &   0  2 13 14  &  0  2 15 17 \\
  0  2 16 19  &  0  2 18 22 &   0  2 20 21 &   0  5  6 20 &   0  5 12 22  &  0  5 13 17\\
  0  5 16 18  &  0  5 19 21 &   0  6  8 19 &   0  6 11 18 &   0  6 12 13  &  0  6 15 21\\
  0  6 16 22  &  0 11 13 22 &   0 12 14 19 &   2  3 11 17 &   2  3 12 22  &  2  3 13 18\\
  2  5  7 19  &  2  5  9 16 &   2  5 11 22 &   2  5 12 14 &   2  5 13 20  &  2  5 15 18\\
  2  5 17 21  &  2  6  7 12 &   2  6 10 21 &   2  6 11 20 &   2  6 14 16  &  2  6 17 22\\
  2  6 18 19  &  2  7 13 17 &   2  7 14 21 &   2  7 15 16 &   2  7 20 22  &  2 11 12 19\\
  2 12 16 18  &  5  6  7 16 &   5  6 12 21 &   5  6 13 18 &   5  6 14 17  &  5  6 19 22\\
  5  7 11 13  &  5  7 14 15 &   5  7 20 21 &   5  8 11 18 &   5  8 12 20  &  5  8 13 19\\
  5 11 14 16  &  5 11 17 19 &   5 12 15 16 &   5 13 16 21 &   5 14 18 22  &  5 15 19 20\\
  5 17 20 22  & 11 12 13 17 &  11 12 18 21 &  11 13 20 21 &  11 14 17 21  & 11 17 18 22\\
 11 19 20 22
\end{tabular}
\end{center}
\qed

The following
lemma was proved by Stern and Lenz in \cite{SL}.

\begin{theorem}
\label{SternLenz} {\rm \cite{SL}} Let $G(L)$ be a graph with vertex set
$Z_{2k}$ where $L$ is a set of integers in the range
$1,2,\ldots,k$, such that $\{a,b\}$ is an edge of $G(L)$ if and only
if $|b-a|\in L$, where $|b-a|=b-a$ if $0\leq b-a\leq k$ and
$|b-a|=a-b$ if $k<b-a<2k$. Then $G(L)$ has a one-factorization if and
only if $2k/gcd (j,2k)$ is even for some $j\in L$.
\end{theorem}

\begin{lemma}
\label{MPQS(35)} There is an MPQS$(35)$.
\end{lemma}

\proof We shall construct an MPQS$(35)$ on $Z_{24}\cup
\{x_1,x_2,\dots,x_{11}\}$.
Beside the blocks of an MPQS$(11)$ on $\{x_1,x_2,\dots,x_{11}\}$, the
other blocks are divided into two parts described below.

For $1\leq i\leq 11$ with $i\neq 8,12$, let $\{F_i, F_{24-i}\}$ be a one-factorization of the graph $G(\{i\})$ over $Z_{24}$, and let $F_{12}$ be the one-factor of the graph $G(\{12\})$ over $Z_{24}$. These one-factorizations exist by Theorem \ref{SternLenz}.

Let $A$ be an $11\times 11$ array as follows.
\begin{center}
\begin{tabular}{rrrrrrrrrrrrrrrrr}

  2 & 12 & 1 & 23 & 3 & 21 & 4 & 20 & 10 & 14 & 22 \\
 12 & 22 & 23 & 1 & 21& 3 & 20 & 4 & 14 & 10 & 2\\
  1 & 23 & 4 & 12 & 2 & 22 & 7 & 17 & 6 & 18 & 20\\
  23 & 1 & 12 & 20 & 22& 2 & 17 & 7 & 18 & 6 & 4\\
  3 & 21 & 2 & 22 & 5 & 12 & 9 & 15 & 1 & 23 & 19\\
  21 & 3 & 22 & 2 & 12 & 19 & 15 & 9 & 23 & 1 & 5\\
   4 & 20 & 7 & 17 & 9 & 15 & 6 & 12 & 2 & 22 & 18\\
   20 & 4 & 17 & 7 & 15 & 9 & 12 & 18 & 22 & 2 & 6\\
  10 & 14 & 6 & 18 & 1 & 23 & 2 & 22 & 9 & 12 & 15\\
  14 & 10 & 18 & 6 & 23 & 1 & 22 & 2 & 12 & 15 & 9\\
  22 & 2 & 20 & 4 & 19 & 5 & 18 & 6 & 15 & 9 & 12\\

\end{tabular}
\end{center}
The first part consists of the following blocks:
 $$\{x_i,x_j,a,b\},\ 1\leq i<j\leq 11,\ \{a,b\}\in F_{A(i,j)}.$$
The blocks in the second part are generated by the following base
blocks modulo $24$.
\begin{center}
\begin{tabular}{llllll}
  $x_1$ 0 5 11 & $x_1$ 0 7 15 & $x_2$ 0 6 11 & $x_2$ 0 8 15 & $x_3$ 0 3 11 & $x_3$ 0 5 14\\
  $x_4$ 0 8 11 & $x_4$ 0 9 14 & $x_5$ 0 4 11 & $x_5$ 0 6 14 & $x_6$ 0 7 11 & $x_6$ 0 8 14\\
  $x_7$ 0 1 11 & $x_7$ 0 3 8  & $x_8$ 0 10 11 & $x_8$ 0 5 8 & $x_9$ 0 3 7 & $x_9$ 0 5 13 \\
  $x_{10}$ 0 4 7 & $x_{10}$ 0 8 13 & $x_{11}$ 0 3 13 & $x_{11}$ 0 1 8 &
  0  1  2  5 & 0  1  3 17 \\ 0  1  6 10 & 0  1  7 18 & 0  1  9 21 &  0  1 13 15 &
  0  1 16 20 & 0  1 19 22 \\ 0  2  4 15 & 0  2  6  8 & 0  2  7  9 & 0  2 10 14 &
  0  3  9 15 & 0  3 14 18 \\ 0  5 10 17 \\

\end{tabular}
\end{center}

It is easy to check that the obtained blocks have no common triples.
So, these blocks form a PQS$(35)$. Further, it has $35+{11\choose
2}\times 12+37\times 24=1583=J(35,4,4)$ blocks and this
PQS$(35)$ is also optimal. Here, we also list the triples
that are not contained in any block so that this construction of
an MPQS$(35)$ is more readable.

\begin{center}
\begin{tabular}{llllll}
$\{x_i,\ a,\ b\}$, & {\rm where} $\{a,b\}\in F_{A(i,i)}$ and
$1\leq i\leq 11$\\
$\{k,\ k+8,\ k+16\}$, & {\rm where} $0\leq k\leq 7$\\
$\{j,\ j+1, j+12\},$ $\{j,\ j+3,\ j+10\}$, & {\rm
where}
$j\in Z_{24}$\\
unused triples of an MPQS$(11)$ & on $\{x_1,x_2,\ldots,x_{11}\}$.
\end{tabular}
\end{center}
\qed

\begin{lemma}
\label{MPQS(47)} There is an MPQS$(47)$.
\end{lemma}

\proof We shall construct an MPQS$(47)$ on $Z_{36}\cup
\{x_1,x_2,\dots,x_{11}\}$.
Beside the blocks of an MPQS$(11)$ on $\{x_1,x_2,\dots,x_{11}\}$, the
other blocks are divided into two parts described below.

For $1\leq i\leq 18$ with $i\neq 4,8,12,16,18$, let $\{F_i, F_{36-i}\}$ be a one-factorization of the graph $G(\{i\})$ over $Z_{36}$, and let $F_{18}$ be the one-factor of the graph $G(\{18\})$ over $Z_{36}$. These one-factorizations exist by Theorem \ref{SternLenz}.

Let $A$ be an $11\times 11$ array as follows.
\begin{center}
\begin{tabular}{rrrrrrrrrrrrrrrrr}
 1 & 18 & 2 & 34 & 3 & 33 & 5 & 31 & 6 & 30 & 35\\
 18 & 35 & 34 & 2 & 33 & 3 & 31 & 5 & 30 & 6 & 1\\
 2 & 34 & 5 & 18 & 1 & 35 & 3 & 33 & 10 & 26 & 31\\
 34 & 2 & 18& 31 & 35 & 1 & 33 & 3 & 26 & 10 & 5\\
 3 & 33 & 1 & 35 & 9 & 18 & 6 & 30 & 2 & 34 & 27\\
 33 & 3 & 35 & 1 & 18 & 27& 39 & 6 & 34 & 2 & 9\\
 5 & 31 & 3 & 33 & 6 & 30 & 10 & 18 & 7 & 29 & 26\\
 31 & 5 & 33 & 3 & 39 & 6 & 18 & 26 & 29 & 7 & 10\\
 6 & 30 & 10 & 26 & 2 & 34 & 7 & 29 & 14 & 18 & 22\\
 30 & 6 & 26 & 10 & 34 & 2 & 29 & 7 & 18 & 22 & 14\\
 35 & 1 & 31 & 5 & 27 & 9 & 26 & 10 & 22 & 14 & 18\\

\end{tabular}
\end{center}
The first part consists of the following blocks:
 $$\{x_i,x_j,a,b\},\ 1\leq i<j\leq 11,\ \{a,b\}\in F_{A(i,j)}.$$
The blocks in the second part are generated by the following base
blocks modulo $36$, where the underlined base block generates 18 distinct blocks.
\begin{center}
\begin{tabular}{llllll}

  $x_1$ 0 4 14 & $x_1$ 0 7 19 & $x_1$ 0 8 21 & $x_1$ 0 9 20 & $x_2$ 0 10 14 & $x_2$ 0 12 19\\
  $x_2$ 0 13 21 & $x_2$ 0 11 20 & $x_3$ 0 4 13 & $x_3$ 0 6 17 & $x_3$ 0 7 21 & $x_3$ 0 8 20\\
  $x_4$ 0 9 13 & $x_4$ 0 11 17  & $x_4$ 0 14 21 & $x_4$ 0 12 20 &  $x_5$ 0 4 17 & $x_5$ 0 5 16\\
  $x_5$ 0 7 15 & $x_5$ 0 10 22  & $x_6$ 0 13 17 & $x_6$ 0 11 16 & $x_6$ 0 8 15 & $x_6$ 0 12 22 \\

   $x_7$ 0 1 9 & $x_7$ 0 2 16  & $x_7$ 0 4 19 & $x_7$ 0 11 23 & $x_8$ 0 8 9 & $x_8$ 0 14 16 \\
   $x_8$ 0 15 19 & $x_8$ 0 12 23  & $x_9$ 0 1 17 & $x_9$ 0 3 12 & $x_9$ 0 4 15 & $x_9$ 0 5 13 \\
  $x_{10}$ 0 16 17 & $x_{10}$ 0 9 12 & $x_{10}$ 0 11 15 & $x_{10}$ 0 8 13 &  $ x_{11}$ 0 3 15 &  $ x_{11}$ 0 2 32 \\
  $ x_{11}$ 0 7 23 & $ x_{11}$ 0 8 25 &  \underline{0  5 18 23} &
  0  1  2 19  &  0  2  4 20 &   0  1  3 29  \\
  0  1  4  5  &  0  1  6 10 &   0  1  7 25 &
  0  1  8 34  &  0  1 11 13 &   0  1 12 15 \\
  0  1 14 31  &  0  1 16 22 &   0  1 21 27 &
  0  1 23 30  &  0  1 24 26 &   0  2  5  7 \\
  0  2  8 30  &  0  2  9 29 &   0  2 11 21 &
  0  2 14 23  &  0  2 15 24 &   0  2 17 27 \\
  0  3  6 31  &  0  3  7 11 &   0  3  8 32 &
  0  3  9 14  &  0  3 13 16 &  0  3 17 20 \\
  0  3 21 30  &  0  4  9 28 &   0  4 10 16 &
  0  4 11 22  &  0  4 18 24 &   0  5 10 17 \\
  0  5 14 22  &  0  5 15 20 &   0  6 13 23 &   0  8 16 26  \\
\end{tabular}
\end{center}

It is easy to check that the obtained blocks have no common triples.
So, these blocks form a PQS$(47)$. Further, it has $35+{11\choose
2}\times 18+81\times 36+18=3959=J(47,4,4)$ blocks and this
PQS$(47)$ is also optimal. Here, we also list the triples
that are not contained in any block so that this construction of
an MPQS$(47)$ is more readable.

\begin{center}
\begin{tabular}{llllll}
$\{x_i, a, b\}$, & {\rm where} $\{a,b\}\in F_{A(i,i)}$ and
$1\leq i\leq 11$,\\
$\{k, k+12, k+24\}$, & {\rm where} $0\leq k\leq 11$,\\
$\{j, j+3, j+18\},$ $\{j, j+2, j+6\}$, & $\{j, j+7, j+20\},$ $\{j, j+8, j+19\}$, {\rm
where} $j\in Z_{36}$,\\
unused triples of an MPQS$(11)$ & on $\{x_1,x_2,\ldots,x_{11}\}$.
\end{tabular}
\end{center}
\qed

 Let $(X,\B)$ be a PQS$(n)$. If there is an $m$-subset $Y$ of $X$ such that every
triple of $Y$ is not contained in any block, then such a PQS is
called a {\it holey PQS} with a {\it hole} $Y$ and denoted by
HPQS$(n,m)$.

\begin{lemma}
\label{MPQS(59)} There is an MPQS$(59)$.
\end{lemma}

\proof We shall construct an MPQS$(59)$ on $Z_{48}\cup
\{x_1,x_2,\dots,x_{11}\}$. The
required blocks are divided into four parts described below.

The first part consists of blocks of an MPQS$(11)$ on $\{x_1,x_2,\dots,x_{11}\}$. For $j\in Z_4$, construct an HPQS$(17,5)$ on $\{4i+j:i\in Z_{12}\}\cup \{x_7,x_8,x_9,x_{10},x_{11}\}$ with $\{x_7,x_8,x_9,x_{10},x_{11}\}$ as a hole and with $J(17,4,4)-J(5,4,4)=156$ blocks. Such a design exists by \cite[Lemma 2.3]{Ji2006}.
The blocks of these four HPQS$(17,5)$ form the second part of blocks.

For $1\leq i\leq 48$ with $i\neq 16,24$, let $\{F_i, F_{48-i}\}$ be a one-factorization of the graph $G(\{i\})$ over $Z_{48}$, and let $F_{24}$ be the one-factor of the graph $G(\{24\})$ over $Z_{48}$. These one-factorizations exist by Theorem \ref{SternLenz}.

Let $A$ be an $11\times 11$ array as follows, where some entries are empty.

\begin{center}
\begin{tabular}{rrrrrrrrrrrrrrrrr}
 3 & 24 & 4 & 44 & 6 & 42 & 1 & 47 & 2 & 46 & 45 \\
 24 & 3 & 44 & 4 & 42 & 6 & 47 & 1 & 46 & 2 & 3\\
 4 & 44 & 5 & 24 & 8 & 40 & 2 & 46 & 1 & 47 & 43\\
 44 & 4 & 24 & 43 & 40 & 8 & 46 & 2 & 47 & 1 & 5\\
 6 & 42 & 8 & 40 & 10 & 24 & 3 & 45 & 5 & 43 & 38\\
 42 & 6 & 40 & 8 & 24 & 38 & 45 & 3 & 43 & 5 & 10\\
 1 & 47 & 2 & 46 & 3 & 45 & & & & & \\
 47 & 1 & 46 & 2 & 45 & 3 & & & & & \\
 2 & 46 & 1 & 47 & 5 & 43 & & & & &\\
 46 & 2 & 47 & 1 & 43 & 5 & & & & & \\
 45 & 3 & 43 & 5 & 38 & 10 & & & & & \\

\end{tabular}
\end{center}
The third part consists of the following blocks:
 $$\{x_i,x_j,a,b\},\ 1\leq i<j\leq 11, (i,j)\not \in \{(i',j'):7\leq i'<j'\leq 11\},\ \{a,b\}\in F_{A(i,j)}.$$
The blocks in the fourth part are generated by the following base
blocks modulo $48$.
\begin{center}
\begin{tabular}{llllll}

  $x_1$ 0 5  12 &    $x_1$ 0  8 22 &    $x_1$ 0   9 27 &    $x_1$ 0  10 25 &    $x_1$ 0  11 28 &    $x_1$ 0  13 29\\
   $x_2$ 0  7 12 &  $x_2$ 0  14 22 &  $x_2$ 0  18 27 &  $x_2$ 0 15 25 &  $x_2$ 0 17 28 &  $x_2$ 0 16 29 \\
   $x_3$ 0  3  9   &  $x_3$ 0   7 26   &  $x_3$ 0 10 28   &  $x_3$ 0 11 23   &  $x_3$ 0 13 27   &  $x_3$ 0  15 31\\
    $x_4$ 0 6 9  &  $x_4$ 0 19 26  &  $x_4$ 0 18 28  &  $x_4$ 0  12 23  &  $x_4$ 0 14 27  &  $x_4$ 0 16 31\\
    $x_5$ 0  1 12  &  $x_5$ 0  2 16  &  $x_5$ 0  4 25  &  $x_5$ 0   7 22  &  $x_5$ 0  9 28  &  $x_5$ 0 13 30\\
    $x_6$ 0  11 12  &  $x_6$ 0  14 16  &  $x_6$ 0 21 25  &  $x_6$ 0  15 22  &  $x_6$ 0 19 28  &  $x_6$ 0  17 30\\
     $x_7$ 0  5  11  &  $x_7$ 0   7 25  &  $x_7$ 0  9 22  &  $x_7$ 0  10 27  &  $x_7$ 0 14 29  &   $x_8$ 0  6 11 \\
   $x_8$ 0 18 25  &  $x_8$ 0 13 22  &  $x_8$ 0 17 27  &  $x_8$ 0 15 29 &    $x_9$ 0  3 21  &  $x_9$ 0   6 19  \\
    $x_9$ 0  7 17  &  $x_9$ 0  9 23  &  $x_9$ 0  11 26 &   $x_{10}$ 0  18 21  &  $x_{10}$ 0  13 19  &  $x_{10}$ 0  10 17  \\
    $x_{10}$ 0  14 23  &  $x_{10}$ 0  15 26 &   $x_{11}$ 0 1 15 & $x_{11}$ 0 2 23 & $x_{11}$ 0  6  13 & $x_{11}$ 0  9 26 \\
     $x_{11}$ 0 11 29  &  0  1  5  6  &  0  1  7  8  &  0  1  9 10 &
  0  1 11 13  &  0  1 14 17  \\  0  1 16 18  &  0  1 19 20 &
  0  1 21 22  &  0  1 23 26  &  0  1 31 33  &  0  1 32 35 \\
  0  1 36 38 &   0  2  5  7  &  0  2  6  8  &  0  2  9 11 &
  0  2 10 15  &  0  2 14 20  \\  0  2 19 21  &  0  2 22 28 &
  0  2 30 36 &   0  2 35 40  &  0  3  7 40 &   0  3  8 37 \\
  0  3 10 39  &  0  3 11 44  &  0  3 12 41 &   0  3 13 38 &
  0  3 14 43  &  0  3 15 18  \\  0  3 19 22  &  0  3 20 23 &
  0  4  9 43  &  0  4 10 14  &  0  4 11 39  &  0  4 13 41 \\
  0  4 17 21  &  0  4 18 22  &  0  4 19 23  &  0  5 16 26 &
  0  5 17 22  &  0  5 18 23  \\  0  5 20 33  &  0  5 21 28 &
  0  5 25 32 &   0  5 27 37  &  0  6 14 37  &  0  6 15 21 \\
  0  6 16 22 &   0  6 17 40  &  0  6 23 29  &  0  7 16 39 &
  0  7 18 34  &  0  7 21 37  \\  0  8 17 29  &  0  8 18 26 &
  0  8 21 33 &   0  8 23 35  &  0  8 27 39  &  0 10 22 36 \\
  0  1  2 25 &   0  2  4 26  &  0  3  6 27  &  0  5 10 29  &
  0  6 12 30 &   0  7 14 31  \\  0  9 18 33  &  0 10 20 34  &
  0 11 22 35 &   0  1  3  4  \\

\end{tabular}
\end{center}

It is easy to check that the above blocks have no common triples.
So, these blocks form a PQS$(59)$. Further, it has $35+4\times 156+[{11\choose
2}-{5\choose 2}]\times 24+130\times 48=7979=J(59,4,4)$ blocks and this
PQS$(59)$ is also optimal. Here, we also list the triples
that are not contained in any block so that this construction of
an MPQS$(59)$ is more readable.

\begin{center}
\begin{tabular}{llllll}
$\{x_i, a, b\}$, & {\rm where} $\{a,b\}\in F_{A(i,i)}$ and
$1\leq i\leq 6$\\

$\{j, j+14, j+15\},$ $\{j, j+21, j+23\}$,& {\rm where} $j\in Z_{48}$,\\
 $\{j, j+7, j+13\},$ $\{j, j+17, j+26\}$,\\
  $\{j, j+18, j+29\}$,\\

unused triples of an MPQS$(11)$ & on $\{x_1,x_2,\ldots,x_{11}\}$,\\
unused triples of four HPQS$(17,5)$ & on $\{4i+j:i\in Z_{12}\}\cup \{x_{7},x_{8},\ldots,x_{11}\}$, $j\in Z_4$.
\end{tabular}
\end{center}
\qed

\begin{lemma}
\label{MPQS(71)} There is an MPQS$(71)$.
\end{lemma}

\proof We shall construct an MPQS$(71)$ on $Z_{48}\cup
\{x_1,x_2,\dots,x_{23}\}$. The
required blocks are divided into four parts described below.

The first part consists of blocks in an MPQS$(23)$ on $\{x_1,x_2,\dots,x_{23}\}$. For $j\in Z_4$, construct an HPQS$(17,5)$ on $\{4i+j:i\in Z_{12}\}\cup \{x_{19},x_{20},x_{21},x_{22},x_{23}\}$ with $\{x_{19},x_{20},x_{21},x_{22},x_{23}\}$ as a hole and with $J(17,4,4)-J(5,4,4)=156$ blocks.  Such a design exists by \cite[Lemma 2.3]{Ji2006}. The blocks of these four HPQS$(17,5)$ form the second part of blocks.

For $1\leq i\leq 48$ with $i\neq 16,24$, let $\{F_i, F_{48-i}\}$ be a one-factorization of the graph $G(\{i\})$ over $Z_{48}$, and let $F_{24}$ be the one-factor of the graph $G(\{24\})$ over $Z_{24}$. These one-factorizations exist by Theorem \ref{SternLenz}.

Let $A$ be a $23\times 23$ array as follows, where some entries are empty.

{\small
\begin{center}
\begin{tabular}{rrrrrrrrrrrrrrrrrrrrrrrrrrrrrrrr}
3 & 24 & 4 & 44 & 6 & 42 & 7 & 41 & 5 & 43 & 8 & 40 & 11 & 37 & 14 & 34 & 23 & 25 & 1 & 47 & 2 & 46 & 45 \\
24 & 45 & 44 & 4 & 42 & 6 & 41 & 7 & 43 & 5 & 40 & 8 & 37 & 11 & 34 & 14 & 25 & 23 & 47 & 1 & 46 & 2 & 3 \\
4 & 44 & 5 & 24 & 7 & 41 & 6 & 42 & 3 & 45 & 9 & 39 & 8 & 40 & 21 & 27 & 18 & 30 & 2 & 46 & 1 & 47 & 43 \\
44 & 4 & 24 & 43 & 41 & 7 & 42 & 6 & 45 & 3 & 39 & 9 & 40 & 8 & 27 & 21 & 30 & 18 & 46 & 2 & 47 & 1 & 5 \\
6 & 42 & 7 & 41 & 1 & 24 & 2 & 46 & 4 & 44 & 10 & 38 & 12 & 36 & 18 & 30 & 20 & 28 & 3 & 45 & 5 & 43 & 47 \\
42 & 6 & 41 & 7 & 24 & 47 & 46 & 2 & 44 & 4 & 38 & 10 & 36 & 12 & 30 & 18 & 28 & 20 & 45 & 3 & 43 & 5 & 1 \\
7 & 41 & 6 & 42 & 2 & 46 & 9 & 24 & 1 & 47 & 4 & 44 & 13 & 35 & 23 & 25 & 15 & 33 & 5 & 43 & 3 & 45 & 39 \\
41 & 7 & 42 & 6 & 46 & 2 & 24 & 39 & 47 & 1 & 44 & 4 & 35 & 13 & 25 & 23 & 33 & 15 & 43 & 5 & 45 & 3 & 9 \\
5 & 43 & 3 & 45 & 4 & 44 & 1 & 47 & 2 & 24 & 11 & 37 & 14 & 34 & 20 & 28 & 19 & 29 & 6 & 42 & 7 & 41 & 46 \\
43 & 5 & 45 & 3 & 44 & 4 & 47 & 1 & 24 & 46 & 37 & 11 & 34 & 14 & 28 & 20 & 29 & 19 & 42 & 6 & 41 & 7 & 2 \\
8 & 40 & 9 & 39 & 10 & 38 & 4 & 44 & 11 & 37 & 13 & 24 & 2 & 46 & 15 & 33 & 17 & 31 & 7 & 41 & 6 & 42 & 35 \\
40 & 8 & 39 & 9 & 38 & 10 & 44 & 4 & 37 & 11 & 24 & 35 & 46 & 2 & 33 & 15 & 31 & 17 & 41 & 7 & 42 & 6 & 13 \\
11 & 37 & 8 & 40 & 12 & 36 & 13 & 35 & 14 & 34 & 2 & 46 & 15 & 24 & 1 & 47 & 5 & 43 & 9 & 39 & 10 & 38 & 33 \\
37 & 11 & 40 & 8 & 36 & 12 & 35 & 13 & 34 & 14 & 46 & 2 & 24 & 33 & 47 & 1 & 43 & 5 & 39 & 9 & 38 & 10 & 15 \\
14 & 34 & 21 & 27 & 18 & 30 & 23 & 25 & 20 & 28 & 15 & 33 & 1 & 47 & 22 & 24 & 2 & 46 & 17 & 31 & 19 & 29 & 26 \\
34 & 14 & 27 & 21 & 30 & 18 & 25 & 23 & 28 & 20 & 33 & 15 & 47 & 1 & 24 & 26 & 46 & 2 & 31 & 17 & 29 & 19 & 22 \\
23 & 25 & 18 & 30 & 20 & 28 & 15 & 33 & 19 & 29 & 17 & 31 & 5 & 43 & 2 & 46 & 14 & 24 & 22 & 26 & 21 & 27 & 34 \\
25 & 23 & 30 & 18 & 28 & 20 & 33 & 15 & 29 & 19 & 31 & 17 & 43 & 5 & 46 & 2 & 24 & 34 & 26 & 22 & 27 & 21 & 14 \\
1 & 47 & 2 & 46 & 3 & 45 & 5 & 43 & 6 & 42 & 7 & 41 & 9 & 39 & 17 & 31 & 22 & 26 &  &  &  &  &  \\
47 & 1 & 46 & 2 & 45 & 3 & 43 & 5 & 42 & 6 & 41 & 7 & 39 & 9 & 31 & 17 & 26 & 22  &  &  &  &  &  \\
2 & 46 & 1 & 47 & 5 & 43 & 3 & 45 & 7 & 41 & 6 & 42 & 10 & 38 & 19 & 29 & 21 & 27  &  &  &  &  & \\
46 & 2 & 47 & 1 & 43 & 5 & 45 & 3 & 41 & 7 & 42 & 6 & 38 & 10 & 29 & 19 & 27 & 21  &  &  &  &  &  \\
45 & 3 & 43 & 5 & 47 & 1 & 39 & 9 & 46 & 2 & 35 & 13 & 33 & 15 & 26 & 22 & 34 & 14  &  &  &  &  &
\end{tabular}
\end{center}
}

The third part consists of the following blocks:
 $$\{x_i,x_j,a,b\},\ 1\leq i<j\leq 23, (i,j)\not \in \{(i',j'):19\leq i'<j'\leq 23\}, \ \{a,b\}\in F_{A(i,j)}.$$
The blocks in the fourth part are generated by the following base
blocks modulo $48$.
\begin{center}
\begin{tabular}{llllll}
 $x_{1}$ 0 9 26 & $x_{1}$ 0 10 28 & $x_{1}$ 0 12 27 & $x_{1}$ 0 13 29
& $x_{2}$ 0 17 26 & $x_{2}$ 0 18 28 \\ $x_{2}$ 0 15 27 & $x_{2}$ 0 16 29
& $x_{3}$ 0 10 26 & $x_{3}$ 0 11 25 & $x_{3}$ 0 12 29 & $x_{3}$ 0 13 28
\\ $x_{4}$ 0 16 26 & $x_{4}$ 0 14 25 & $x_{4}$ 0 17 29 & $x_{4}$ 0 15 28
& $x_{5}$ 0 8 25 & $x_{5}$ 0 9 22 \\ $x_{5}$ 0 11 27 & $x_{5}$ 0 14 29
& $x_{6}$ 0 17 25 & $x_{6}$ 0 13 22 & $x_{6}$ 0 16 27 & $x_{6}$ 0 15 29
\\ $x_{7}$ 0 8 27 & $x_{7}$ 0 10 22 & $x_{7}$ 0 11 28 & $x_{7}$ 0 14 30
& $x_{8}$ 0 19 27 & $x_{8}$ 0 12 22 \\ $x_{8}$ 0 17 28 & $x_{8}$ 0 16 30
& $x_{9}$ 0 8 26 & $x_{9}$ 0 9 21 & $x_{9}$ 0 10 23 & $x_{9}$ 0 15 31
\\ $x_{10}$ 0 18 26 & $x_{10}$ 0 12 21 & $x_{10}$ 0 13 23 & $x_{10}$ 0 16 31
& $x_{11}$ 0 1 19 & $x_{11}$ 0 3 23 \\ $x_{11}$ 0 5 21 & $x_{11}$ 0 12 26
& $x_{12}$ 0 18 19 & $x_{12}$ 0 20 23 & $x_{12}$ 0 16 21 & $x_{12}$ 0 14 26
\\ $x_{13}$ 0 3 19 & $x_{13}$ 0 4 21 & $x_{13}$ 0 6 26 & $x_{13}$ 0 7 25
& $x_{14}$ 0 16 19 & $x_{14}$ 0 17 21 \\ $x_{14}$ 0 20 26 & $x_{14}$ 0 18 25
& $x_{15}$ 0 3 12 & $x_{15}$ 0 4 11 & $x_{15}$ 0 5 13 & $x_{15}$ 0 6 16
\\ $x_{16}$ 0 9 12 & $x_{16}$ 0 7 11 & $x_{16}$ 0 8 13 & $x_{16}$ 0 10 16
& $x_{17}$ 0 1 13 & $x_{17}$ 0 3 11 \\ $x_{17}$ 0 4 10 & $x_{17}$ 0 7 16
& $x_{18}$ 0 12 13 & $x_{18}$ 0 8 11 & $x_{18}$ 0 6 10 & $x_{18}$ 0 9 16
\\ $x_{19}$ 0 10 25 & $x_{19}$ 0 11 29 & $x_{19}$ 0 13 27
& $x_{20}$ 0 15 25 & $x_{20}$ 0 18 29 & $x_{20}$ 0 14 27
\\ $x_{21}$ 0 9 23 & $x_{21}$ 0 11 26 & $x_{21}$ 0 13 30
& $x_{22}$ 0 14 23 & $x_{22}$ 0 15 26 & $x_{22}$ 0 17 30
\\ $x_{23}$ 0 6 25 & $x_{23}$ 0 7 18 & $x_{23}$ 0 10 27
 & 0  1  2 25 &   0  2  4 26  &  0  3  6 27  \\  0  5 10 29  &
  0  6 12 30 &   0  7 14 31  &  0  9 18 33  &  0 10 20 34  &
  0 11 22 35 \\
  0  1  3  4 &   0  1  5  6 &   0  1  7  8 &   0  1  9 10 &
  0  1 11 12 &   0  1 14 15 \\   0  1 16 17 &   0  1 18 20 &
  0  1 21 22 &   0  1 23 26 &   0  1 29 31 &   0  2  5  7 \\
  0  2  6  8 &   0  2  9 11 &   0  2 10 37 &   0  2 12 14 &
  0  2 13 40 &   0  2 15 17 \\   0  2 16 18 &   0  2 21 28 &
  0  2 22 29 &   0  2 23 27 &   0  3  7 21 &   0  3  8 38 \\
  0  3  9 20 &   0  3 10 41 &   0  3 13 43 &   0  3 14 17 &
  0  3 15 18 &   0  3 16 22 \\   0  3 29 35 &   0  3 30 44 &
  0  3 31 42 &   0  4  9 30 &   0  4 13 17 &   0  4 14 19 \\
  0  4 15 23 &   0  4 22 43 &   0  4 29 37 &   0  4 33 38 &
  0  5 11 16 &   0  5 12 41 \\   0  5 14 22 &   0  5 17 23 &
  0  5 20 25 &   0  5 30 36 &   0  5 31 39 &   0  6 13 41 \\
  0  6 14 20 &   0  6 15 21 &   0  7 15 22 &   0  7 23 32 &
  0  9 19 28 &   0 11 23 36 \\

\end{tabular}
\end{center}

It is easy to check that the above blocks have no common triples.
So, these blocks form a PQS$(71)$. Further, it has $419+4\times 156+[{23\choose
2}-{5\choose 2}]\times 24+150\times 48=14075=J(71,4,4)$ blocks and this
PQS$(71)$ is also optimal. Here, we also list the triples
that are not contained in any block so that this construction of
an MPQS$(71)$ is more readable.

\begin{center}
\begin{tabular}{llllll}
$\{x_i, a, b\}$, & {\rm where} $\{a,b\}\in F_{A(i,i)}$ and
$1\leq i\leq 18$,\\
$\{j, j+19, j+25\},$ $\{j, j+11, j+18\}$, & $\{j, j+17, j+27\},$ {\rm
where} $j\in Z_{48}$,\\
unused triples of an MPQS$(23)$ & on $\{x_{1},x_{2},\ldots,x_{23}\}$, \\
unused triples of four HPQS$(17,5)$ & on $\{4i+j:i\in Z_{12}\}\cup \{x_{19},x_{20},\ldots,x_{23}\}$, $j\in Z_4$.
\end{tabular}
\end{center}
\qed

\section{Constructions for MPQSs}

In this section we describe recursive constructions for MPQS$(n)$'s
via candelabra quadruple systems.

Let $v$ be a non-negative integer, let $t$ be a positive integer
and let $K$ be a set of positive integers. A {\it candelabra
$t$-system} (or $t$-$CS$) of order $v$, and
block sizes from $K$ is a quadruple $(X, S, \G, \A)$ that
satisfies the following properties:

 \medskip

(1) $X$ is a set of $v$ elements (called {\it points}).

(2) $S$ is a subset (called the {\it stem} of the candelabra) of
$X$ of size $s$.

(3) $\G=\{G_1,G_2,\ldots\}$ is a set of non-empty subsets
 (called {\it groups} or {\it branches}) of $X\backslash S$,
 which partition $X\backslash S$.

(4) $\A$ is a family of subsets (called {\it blocks}) of $X$,
 each of cardinality from $K$.

(5) Every $t$-subset $T$ of $X$ with $|T\cap (S\cup G_i)|<t$ for
all $i$,
 is contained in a unique block and no $t$-subsets of $S\cup G_i$
 for all $i$, are contained in any block.

 \medskip

\noindent Such a system is denoted by $CS(t,K,v)$. By the {\it
group type} (or {\it type}) of a $t$-$CS$ $(X,S,\Gamma,\A)$ we
mean the list $(|G||G\in \Gamma: |S|)$ of group sizes and stem
size. The stem size is separated from the group sizes by a colon.
If a $t$-$CS$ has $n_i$ groups of size $g_i$, $1\leq i\leq r$, and
stem size $s$, then we use the notation $(g_1^{n_1}g_2^{n_2}\cdots
g_r^{n_r}: s)$ to denote group type.  A candelabra system with
$t=3$ and $K=\{4\}$ is called a {\it candelabra quadruple system}
and briefly denoted by CQS$(g_1^{n_1}g_2^{n_2}\cdots g_r^{n_r}:
s)$. A $CS(t,K,v)$ with group type $(1^v:0)$ is usually called a
{\it $t$-wise balanced design} and shortly denoted by S$(t,K,v)$.
As well,  the group set ${\cal G}$ and the stem $S$ in the
quadruple $(X,S,{\cal G},{\cal A})$ can be omitted and we write
$(X,{\cal A})$ instead of $(X,S,{\cal G},{\cal A})$. When $K=\{k\}$, we simply write $k$ instead of $K$.

\begin{theorem}{\rm \cite{Mills1974}}
\label{CQS(6^k:0)} There is a CQS$(6^k:0)$ for any $k\geq 0$.
\end{theorem}

\begin{theorem}
\label{CQS(g^3:s)} {\rm \cite{Hartman1982,Hartman1980,Lenz1985}} A
CQS$(g^3:s)$ exits for all even $s$ and all $g\equiv 0,s\ (mod\
6)$ with $g\geq s$.
\end{theorem}

\begin{theorem}
\label{CQS(g^4:s)} {\rm \cite{GH1991,Zhang2011}}
There exists a CQS$(g^4 : s)$ if and only if $g\equiv 0 \pmod 2$, $s\equiv 0 \pmod 2$ and
$0\leq s\leq 2g$.
\end{theorem}

\begin{theorem}
\label{CQS(g^5:s)} {\rm \cite{Zhang2011}}
A CQS$(g^5:s)$ exists for all $g\equiv 0 \pmod 6$, $s\equiv 0 \pmod 2$ and $0\leq s\leq 3g$.
\end{theorem}

\begin{lemma}
\label{CQS(12^k:6)} {\rm \cite{Ji2006}} There is a CQS$(12^k:6)$ for any $k\geq 3$.
\end{lemma}

With the aid of CQSs,  a construction of MPQS$(n)$ for $n\equiv 5\ (mod\ 6)$  has been stated in \cite{Ji2006}.

\begin{construction}
\label{3.1} {\rm \cite{Ji2006}} Suppose that there is a
CQS$(g_0^{1}g_1^{a_1}g_2^{a_2}\cdots g_r^{a_r}:s)$, where $s\equiv
6 \pmod {12}$, $g_i\equiv 0 \pmod {12}$ for $1\leq i\leq r$, and
$g_0\equiv 0 \pmod {6}$. If there is an MPQS$(g_0+s-1)$ and an
HPQS$(g_i+s-1,s-1)$ with $J(g_i+s-1,4,4)-J(s-1,4,4)$ blocks for
$1\leq i\leq r$, then there is an MPQS$(\sum_{1\leq i\leq
r}{a_ig_i}+g_0+s-1)$.
\end{construction}

Similar to the proof of Construction \ref{3.1}, we can get another construction for $n\equiv 5\pmod 6$.

\begin{construction}
\label{3.2} Suppose that there is a
CQS$(g_0^{1}g_1^{a_1}g_2^{a_2}\cdots g_r^{a_r}:s)$, where $s\equiv g_i\equiv 0\pmod {12}$ for $0\leq i\leq r$.
If there is an MPQS$(g_0+s-1)$ and an
HPQS$(g_i+s-1,s-1)$ with $J(g_i+s-1,4,4)-J(s-1,4,4)$ blocks for
$1\leq i\leq r$, then there is an MPQS$(\sum_{1\leq i\leq
r}{a_ig_i}+g_0+s-1)$.
\end{construction}

\proof Let $(X,S,\G,\B)$ be a given CQS$(g_0^{1}g_1^{a_1}g_2^{a_2}\cdots
g_r^{a_r}:s)$. We shall construct the desired design as follows.

Take a point $x$ from $S$ and let $S'=S\setminus \{x\}$. Denote
$\B'=\{B\in \B: x\not\in B\}$. For a special group $G$ with
$|G|=g_0$, construct an MPQS$(g_0+s-1)$ on $G\cup S'$. Such a
design exists by assumption. Denote its block set by $\C_G$. For
each group $G'\neq G$, construct an HPQS$(|G'|+s-1,s-1)$  on
$G'\cup S'$ with a hole $S'$ and $J(|G'|+s-1,4,4)-J(s-1,4,4)$
blocks. Such a design exists by assumption. Denote its block set
by $\C_{G'}$.

Let $$\A=\B'\bigcup \C_G\bigcup (\bigcup_{G'\in \G,G'\neq
G}\C_{G'}).$$

It is easy to see that all blocks in $\A$ have no common triples.
So, $(X\setminus \{x\},\A)$ is a PQS$(\sum_{1\leq i\leq
r}{a_ig_i}+g_0+s-1)$. It is left to check that $|\A|=J(\sum_{1\leq
i\leq r}{a_ig_i}+g_0+s-1,4,4)$.

Let $u=g_0+\sum_{1\leq i\leq r}a_ig_i$ and $\B_x=\{B\in \B: x\in B
\}$. Clearly, $\B'=\B\setminus \B_x$. Since $\B$ is the block set
of a CQS$(g_0^{1}g_1^{a_1}g_2^{a_2}\cdots g_r^{a_r}:s)$ and
$\{B\setminus \{x\}: B\in \B_x\}$ is the block set of a
GDD$(2,3,u)$ of type $g_0^{1}g_1^{a_1}g_2^{a_2}\cdots g_r^{a_r}$,
we have that $|\B|=\frac{1}{4}[{u+s\choose 3}-{g_0+s\choose
3}-\sum_{1\leq i\leq r}a_i({g_i+s\choose 3}-{s\choose 3})]$ and
$|\B_x|=\frac{1}{3}[{u\choose 2}-{g_0\choose 2}-\sum_{1\leq i\leq
r}a_i{g_i\choose 2}]$. By simple computing, we have
$$|\B'|=|\B|-|\B_x|=\frac{1}{24}[u^3-g_0^3-\sum_{1\leq i\leq r}a_ig_i^3+(3s-7)(u^2-g_0^2-\sum_{1\leq i\leq r}a_ig_i^2)].$$
By the definition, $J(n,4,4)=\frac{1}{24}[n^3-4n^2+n-18]$ for $n\equiv 11\pmod {12}$.
Since $|\C_{G'}|=J(|G'|+s-1,4,4)-J(s-1,4,4)$, $|G'|\equiv 0\ (mod\
12)$ and $s-1\equiv 11\pmod {12}$, we have
$|\C_{G'}|=\frac{1}{24}[|G'|^3+|G'|^2(3s-7)+|G'|(3s^2-14s+12)]$.
So, $$|\bigcup_{G'\in \G,G'\neq G}\C_{G'}|=\frac{1}{24}\sum_{1\leq
i\leq r}a_i[g_i^3+g_i^2(3s-7)+g_i(3s^2-14s+12)].$$ Also,
$$|\C_G|=\frac{1}{24}[g_0^3+g_0^2(3s-7)+g_0(3s^2-14s+12)+s^3-7s^2+12s-24].$$

Since $|\A|=|\B'|+|\C_G|+|\bigcup_{G'\in \G,G'\neq G}\C_{G'}|$,
the number of blocks is
$$
\frac{1}{24}[u^3+u^2(3s-7)+u(3s^2-14s+12)+s^3-7s^2+12s-24],
$$
\noindent which is equal to $J(u+s-1,4,4)$. This completes the
proof. \qed

From Constructions \ref{3.1}-\ref{3.2}  CQSs are useful in the constructions
for MPQSs. A recursive construction
for CQSs has been stated in  \cite{Ji2006}.

Let $v$ be a non-negative integer, let $t$ be a positive integer
and $K$ be a set of positive integers. A {\it group divisible
$t$-design} (or $t$-GDD) of order $v$ and block sizes from $K$
denoted by GDD($t,K,v$) is a triple $(X,\G,\B)$ such that

   (1) $X$ is a set of $v$ elements (called {\it points});

   (2) $\G=\{G_1,G_2,\ldots\}$ is a set of non-empty subsets (called {\it groups}) of $X$,
   which  partition $X$;

   (3) $\B$ is a family of subsets (called {\it blocks}) of $X$ each of cardinality
from $K$ such that each block intersects any given group in at
most one point;

   (4) each $t$-set of points from $t$ distinct groups is contained in exactly
   one block.\\
\noindent The {\it type} of $t$-GDD is defined as the list
$\{|G|:G\in \G\}$. When $K=\{k\}$, we simply write $k$ for $K$.

A GDD$(3, 4, v)$ of type $r^m$ is called an {\it H design} (as in
\cite{M}) and denoted by H$(m,r,4,3)$.

\begin{theorem}{\rm \cite {Ji2009,M}}
\label{H-design} For $m>3$ and $m\not= 5$, an H$(m,r,4,3)$ exists if
and only if $rm$ is even and $r(m-1)(m-2)$ is divisible by $3$.
For $m=5$, H$(5,r,4,3)$ exists if $r$ is even, $r\neq 2$ and
$r\not\equiv 10,26 \pmod{48}$.
\end{theorem}

Let $(X,S,\G,\A)$ be a $CS(3,K,v)$  of type
$(g_1^{a_1}g_2^{a_2}\cdots g_r^{a_r}:s)$ with  $s>0$ and let
$S=\{\infty_1,\ldots,\infty_s\}$. For $1\leq i\leq s$, let
$\A_i=\{A\setminus \{\infty_i\}: A\in \A, \infty_i\in A\}$ and
$\A_T=\{A\in \A: A\cap S=\emptyset\}$. Then the ($s+3$)-tuple ($X,
\G, \A_1, \A_2, \ldots, \A_s, \A_T$) is called an {\it $s$-fan
design} (as in \cite{Hartman1994}). If block sizes of $\A_i$ and
$\A_T$ are from $K_i$($1\leq i\leq s$) and $K_T$, respectively,
then the $s$-fan design is denoted by $s$-FG($3, (K_1, K_2,
\ldots, K_s, K_T),\sum_{i=1}^{r}{a_ig_i}$) of type
$g_1^{a_1}g_2^{a_2}\cdots g_r^{a_r}$.

Below is a recursive construction for CQSs, which was obtained by applying Hartman's fundamental construction for $3$-CSs \cite {Hartman1994},
\begin{lemma}
\label{C-CQS} {\rm \cite{Ji2006}} Suppose there is an $e$-FG$(3, (K_1, \cdots,
K_e,K_T), v)$ of type $g_1^{a_1}g_2^{a_2}$ $\cdots g_r^{a_r}$ with
$e\geq 1$, $K_i\subset \{k\geq 3:$ $k$ is an integer$\}$ $(2\leq
i\leq e)$ and $K_T\subset \{k\geq 4:$ $k$ is an integer$\}$.
Suppose that $b\equiv 0\ (mod\ 6 )$ and there exists a
CQS$(b^{k_1}:s)$ for any $k_1\in K_1$. Then there exists a
CQS$({(bg_1)}^{a_1}(bg_2)^{a_2}\cdots (bg_r)^{a_r}:b(e-1)+s)$.
\end{lemma}

In the next section, we shall obtain some CQSs and then
determine the packing numbers $A(n,4,4)$.

\section{Existence of MPQSs}

In this section we shall determine the existence of the last 21 undecide MPQS$(n)$ for
$n\in  \{6k+5:k=3,5,7,9,11,13,15,19,23,25,$ $27,29,31,33,35,45,47,75,77,79,159\}$.

\begin{lemma}
\label{CQS(24^k:12)} There is a CQS$(24^k:12)$ for any $k\geq 3$.
\end{lemma}

\proof For $k\equiv 0,1\pmod  3$, there is a
2-FG$(3,(3,3,4),2k)$ of type $2^{k}$, which can be obtained by
deleting two points from an SQS$(2k+2)$ in \cite {Hanani1960}.
Applying Lemma \ref{C-CQS} with $b=12$ and the known CQS$(12^3:0)$ in
Lemma \ref{CQS(g^3:s)} gives a CQS$(24^k:12)$.

For $k\equiv 2\pmod  3$, there is a
2-FG$(3,(\{3,5\},\{3,5\},\{4,6\}),2k)$ of type $2^k$, which can be
obtained by deleting two points from two distinct groups of a
CQS$(6^{(k+1)/3}:0)$ in Theorem \ref{CQS(6^k:0)}. A
CQS$(24^k:12)$ is then obtained by applying Lemma \ref{C-CQS}
with $b=12$ and the known CQS$(12^j:0)$ ($j=3,5$) by Theorem \ref{CQS(g^3:s)} and Theorem \ref{CQS(g^5:s)}. \qed

\begin{lemma}
\label{MPQS(24k+11)} There is an MPQS$(24k+11)$ for any $k\geq 3$. So, there is an MPQS$(n)$ for $n\in \{6k+5:k=5,13,25,29,33,45,77\}$
\end{lemma}

\proof By Lemma \ref{CQS(24^k:12)}, there is a CQS$(24^k:12)$. Apply
Construction \ref{3.2} with $g_0=g_1=24$, $r=1$, $a_1=k-1$ and
$s=12$. Since there is an
MPQS$(35)$ and an HPQS$(35,11)$ with
$J(35,4,4)-J(11,4,4)$ blocks which exists from the proof of Lemma \ref{MPQS(35)}, there is an MPQS$(24k+11)$. \qed

\begin{lemma}
\label{MPQS(6^27+5)} There is an MPQS$(6k+5)$ for $k\in \{27,35\}$.
\end{lemma}

\proof For $k=27$, there is a CQS$(48^3:24)$ by Theorem \ref{CQS(g^3:s)}. Since there is an MPQS$(71)$ and an HPQS$(71,23)$ with $J(71,4,4)-J(23,4,4)$ blocks which exists from the proof of Lemma \ref{MPQS(71)}, there is an MPQS$(6k+5)$  by Construction \ref{3.2}.

 For $k=35$, there is a CQS$(48^4:24)$ by Theorem \ref{CQS(g^4:s)}. Since there is an MPQS$(71)$ and an HPQS$(71,23)$ with $J(71,4,4)-J(23,4,4)$ blocks, there is an MPQS$(6k+5)$  by Construction \ref{3.2}. \qed

\begin{lemma}
\label{MPQS(6^31+5)} There is an MPQS$(191)$.
\end{lemma}

\proof Deleting one point from an SQS$(16)$ containing a subdesign S$(2,4,16)$ \cite[Theorem 1.3]{Ji2012} gives a 1-FG$(3,(3,4),15)$ of type $3^5$.
Applying Lemma \ref{C-CQS} with $b=12$ and the known CQS$(12^3:12)$
 gives a CQS$(36^5:12)$. Since there is an
MPQS$(47)$ and an HPQS$(47,11)$ with $J(47,4,4)-J(11,4,4)$ blocks which exists from the proof of Lemma \ref{MPQS(47)}, there is an MPQS$(191)$ by Construction \ref{3.2}.\qed

The next lemma is the well-known result on S($3,k,v$)s.

\begin{lemma}
\label{inversive plane} {\rm \cite{Hanani1979}} For any prime
power $q$ there exists an S$(3, q+1, q^2+1)$ and an S$(3,6,22)$.
\end{lemma}

\begin{lemma}
\label{MPQS(6^19+5)} There is an MPQS$(6k+5)$ for $k\in \{19,23\}$.
\end{lemma}

\proof Deleting two points of an S$(3,6,k+3)$ by Lemma \ref{inversive plane} gives a 2-FG$(3,(5,5,6),k+1)$ of type $4^{(k+1)/4}$. Further, deleting one point from a group give a 2-FG$(3,(\{4,5\},\{4,5\},\{4,5,6\}),k)$ of type $4^{(k-3)/4}3^1$.  Applying Lemma \ref{C-CQS} with $b=6$ and the known CQS$(6^j:0)$ for $j\in \{4,5\}$ in
Theorem \ref{CQS(6^k:0)} gives a CQS$(24^{(k-3)/4}18^1:6)$. Since there is an HPQS$(29,5)$ with $J(29,4,4)-J(5,4,4)$ blocks \cite[Lemma 4.4]{Ji2006} and an
MPQS$(23)$ by Lemma \ref{MPQS(23)}, there is an MPQS$(6k+5)$ by Construction \ref{3.1}. \qed

\begin{lemma}
\label{MPQS(6^15+5)} There is an MPQS$(95)$.
\end{lemma}

\proof Deleting one point of an S$(3,5,17)$  by Lemma \ref{inversive plane} gives a 1-FG$(3,(4,5),16)$ of type $4^4$. Further, deleting one point from a group give a 1-FG$(3,(\{3,4\},\{4,5\}),15)$ of type $4^33^1$.  Applying Lemma \ref{C-CQS} with $b=6$ and the known CQS$(6^j:6)$
for $j\in \{3,4\}$  in Theorem \ref{CQS(g^3:s)} and Theorem \ref{CQS(g^4:s)} gives a CQS$(24^318^1:6)$. Since there is an HPQS$(29,5)$ with $J(29,4,4)-J(5,4,4)$ blocks \cite[Lemma 4.4]{Ji2006} and an
MPQS$(23)$ by Lemma \ref{MPQS(23)}, there is an MPQS$(95)$  by Construction \ref{3.1}. \qed

\begin{lemma}
\label{MPQS(6^47+5)} There is an MPQS$(6k+5)$ for $k\in \{47,75,79,159\}$.
\end{lemma}

\proof For $k=47$, deleting two points from an S$(3,8,50)$   by Lemma \ref{inversive plane} gives a 2-FG$(3,(7,7,8),48)$ of type $6^8$. Further, deleting one point gives a 2-FG$(3,(\{6,7\},\{6,7\},\{7,8\}),48)$ of type $6^75^1$.  Applying Lemma \ref{C-CQS} with $b=6$ and the known CQS$(6^j:0)$ for $j\in \{6,7\}$  by Theorem \ref{CQS(6^k:0)}
 gives a CQS$(36^730^1:6)$. Since there is an HPQS$(41,5)$ with $J(41,4,4)-J(5,4,4)$ blocks by  \cite[Lemma 4.4]{Ji2006} and an
MPQS$(35)$ by Lemma \ref{MPQS(35)}, there is an MPQS$(6k+5)$  by Construction \ref{3.1}.

For $k=75,79$, deleting two points from an S$(3,10,82)$   by Lemma \ref{inversive plane} gives a 2-FG$(3,(9,9,10),80)$ of type $8^{10}$. Further, deleting $80-k$ points from one group gives a 2-FG$(3,(\{8,9\},\{8,9\},\{8,9,10\}),75)$ of type $8^9(k-72)^1$.  Applying Lemma \ref{C-CQS} with $b=6$ and the known CQS$(6^j:0)$ for $j\in \{8,9\}$  by Theorem \ref{CQS(6^k:0)}
 gives a CQS$(48^9(6k-432)^1:6)$. Since there is an HPQS$(53,5)$ with $J(53,4,4)-J(5,4,4)$ blocks by \cite[Lemma 4.4]{Ji2006} and an
MPQS$(23)$ by Lemma \ref{MPQS(23)} and an MPQS$(47)$ by Lemma \ref{MPQS(47)}, there is an MPQS$(6k+5)$  by Construction \ref{3.1}.

For $k=159$, deleting two points from an S$(3,14,169)$   by Lemma \ref{inversive plane} gives a 2-FG$(3,(13,13,$ $14),168)$ of type $12^{14}$. Further, deleting nine points from one group gives a 2-FG$(3,(\{12,13\},$ $\{12,13\},\{12,13,14\}),159)$ of type $12^{13}3^1$.  Applying Lemma \ref{C-CQS} with $b=6$ and the known CQS$(6^j:0)$ for $j\in \{12,13\}$  by Theorem \ref{CQS(6^k:0)}
 gives a CQS$(72^{13}18^1:6)$. Since there is an HPQS$(77,5)$ with $J(77,4,4)-J(5,4,4)$ blocks by  \cite[Lemma 4.4]{Ji2006} and an
MPQS$(23)$ by Lemma \ref{MPQS(23)}, there is an MPQS$(6k+5)$  by Construction \ref{3.1}. \qed

Combining Theorem \ref{1.1}, Lemmas \ref{MPQS(23)}-\ref{MPQS(71)}, Lemmas \ref{MPQS(24k+11)}-\ref{MPQS(6^31+5)} and Lemmas \ref{MPQS(6^19+5)}-\ref{MPQS(6^47+5)}, we obtain the main result of this paper.

\begin{theorem}
For any positive integer $n$, it holds that $A(n,4,4)=J(n,4,4)$.
\end{theorem}

\bigskip


\end{document}